\documentclass[10pt]{article}
\usepackage{mathtext}
\usepackage[T2A]{fontenc}
\usepackage[russian,english]{babel}

\usepackage{amssymb}
\usepackage{graphicx}
\begin{document}

\centerline{\bf I.Kh. Sabitov}
\vskip 0.2cm
\begin{center}
{\bf The volume of an infinitesimally flexible polyhedron is \\
a multiple root of its volume polynomial}
\end{center}
\vskip 0.4cm

%{\footnotesize
%Рассматривается класс многогранников, который мы называем пирамидами, и  доказывается, что при некоторых простых, но
%довольно общих условиях на их внешнее строение пирамиды являются неизгибаемыми, причем это свойство неизгибаемости при
%соответствующих предположениях можно установить  и в многомерных пространствах постоянной кривизны.}
%\end{document}

\vskip 0.4cm
{\bf 1.} In [1] among the many open problems there is an one under the number 6 in which the following property of
polyhedra is announced as a suggestion:

\textbf{Theorem.}\emph{Algebraic volume of an infinitesimally non-rigid polyhedron is a multiple root of its volume
 polynomial}.

 Here we present a sketch of a proof of this assertion in three-dimensional space. Because any flexible polyhedron is infinitesimally
 (inf.) bendable one then as a corollary of our theorem we can affirm that the volume of any flexible polyhedron is a multiple
 root of its volume polynomial

{\bf 2.} Let $P$ be a simplicial orientable inf. bendable polyhedron in $R^3$ of any combinatorial
structure $K$. Let points $M_i (x_i, y_i, z_i), 1\le i\le n$, be its vertices. Inf. bendability (i.e., flexibility,
non-rigidity) of the polyhedron $P$ means that there are vectors ${\bf Z}_i=\{\xi_i, \eta_i, \zeta_i\}$ attached to
vertices $M_i$ satisfying the equations
\begin{equation}
\label{1}
(x_i-x_j)(\xi_i-\xi_j)+(y_i-y_j)(\eta_i-\eta_j)+(z_i-z_j)(\zeta_i-\zeta_j=0,
\end{equation}
written for all edges $(i,j)\in E$ of $P$ (where by $E$ we note the set of all edges of $P$, and indices
correspond to numbers $i$ и $j$ of end vertices of edges). Under the deformation
\begin{equation}
\label{2}
(x_i,y_i, z_i)\to (x_i+\varepsilon \xi_i, y_i+\varepsilon \eta_i, z_i+\varepsilon \zeta_i)
\end{equation}
the length $l_{ij}$ of the edge $(i,j)$ changes to an inf. small $o(\varepsilon),\varepsilon\to 0$.

 The  vectors $Z_i$ are found as a solution of homogenous linear system of $|E|$ equations for $3n$ unknowns
 (in reality for $3n-6$ unknowns because 6 unknowns can be fixed by a motion of $P$ as a solid).
 The number of edges is $|E|=3n+6g-6$ where $g$ is topological genus of $P$. So for the existence of a nontrivial
 solution of our system is it necessary and sufficient that the rank of the main matrix of system be less then the number
 of unknowns equal to $3n-6$. So in the case of inf. bendability of $P$ the coordinates
 of its vertices should satisfy to one or more polynomial equations of the form
  \begin{equation}
 \label{3}
 det D(x)=0~(при~g=0), ~ det D_1(x)=0, ..., det D_s(x)=0,~ s=C_{3n+6g-6}^{3n-6}.
 \end{equation}
% \end{document}

 It is important to remark that if an inf. small deformation is \emph{trivial} that is it an initial velocity vector of a
 motion of $P$ as a solid body then the distances between \emph{all} vertices of $P$ are changed by the order
 $o(\varepsilon), \varepsilon\to 0$. With an additional condition the inverse statement is true too:

 \textbf{Lemma 1.}\emph{ If a polyhedron is not a situated on a plane then for any nontrivial inf. deformation there exists
 a small of nonzero length diagonal whose length changes to the exact order $O(\varepsilon), \varepsilon\to 0$.}

 \textbf{Remarks.} 1) A small diagonal is the one between two vertices of two faces with a common edge. In some cases
 this diagonal is in reality an edge. The lemma is equivalent to the affirmation that there is a  dihedral angle
 changing as an exact $O(\varepsilon), \varepsilon\to 0$.

 2) If all vertices of $P$ are situated on a plane  then can  be that under a nontrivial inf. bending the lengths of all
 diagonals are changing to the order $o(\varepsilon), \varepsilon\to 0$.

{\bf 3.} Let's continue the proof. Recall that a volume polynomial of a polyhedron $P$ is any polynomial of the form
 \begin{equation}
 \label{4}
 Q(V, l)=V^{2N}+\sum_1^N a_i(l)V^{2N-2i},
 \end{equation}
 where coefficients $a_i(l)$  are some polynomials too in the set $l$ of squares of lengths of edges of $P$
 such that after the substitution in (\ref{4}) instead the volume $V$ and the squares of lengths of edges
 their expressions in coordinates of vertices $x=(x_1, y_1, z_1,..., z_n)$ the value $Q(V(x),l(x))$ becomes
 identically zero relatively to all coordinates.

 In [1] one can find a proof of existence theorem for such a polynomial with a detailed description of the background history.

 Let's consider a new polyhedron $P_{\varepsilon}$ with vertices coordinates $(x_i+\varepsilon \xi_i.
 y_i+\varepsilon \eta_i, z_i+\varepsilon \zeta_i)$ and with the same combinatorial structure $K$.
 For the squares of lengths of its edges we have $l_{ij}^2(\varepsilon)= l_{ij}^2+\varepsilon^2L_{ij}^2$,
where $L_{ij}^2 = (\xi_i-\xi_j)^2+(\eta_i-\eta_j)^2+ (\zeta_i-\zeta_j)^2$. Renumber all edges by  the index
$k, 1\le k\le |E|$ and compose for $P_{\varepsilon}$ its volume polynomial:
 \begin{equation}
 \label{5}
 Q(V, l,\varepsilon)= V^{2N}(\varepsilon)+\sum_{i=1}^N a_i(l_{\varepsilon})V^{2N-2i}(\varepsilon)=0, \forall \varepsilon.
 \end{equation}
 %\end{document}
 Evidently for $\varepsilon =0$ this polynomial becomes a volume polynomial for the initial polyhedron with $V_0=V(0)$.
 The derivation of (5) by $\varepsilon$ gives
 $$
 Q'_VV'_{\varepsilon}+2\varepsilon\sum_{i=1}^N (\sum_{j=1}^{|E|}\frac{\partial a_i(l_{\varepsilon})}{\partial l_k}L_k^2)
 V^{2N-2i}=0.
 $$
 If $\lim V'_{\varepsilon}\ne 0$ when $\varepsilon \to 0$, then one has $Q'_V(V_0)=0$, so the multiplicity
 of the root $V=V_0$ is proven.

Let now be $V'_{\varepsilon}(0)=0$. Algebraic volume of a polyhedron is defined as the sum of oriented volumes of tetrahedra
 with a common vertex and the bases on oriented faces of the polyhedron. Let this common vertex be taken as a vertex of $P$
 and choose this vertex as the origin for coordinate system. Then the volume of a tetrahedron  with vertices
 $M_i, M_j, M_k$ is given by the formula
 $$
 V_{ijk}(\varepsilon)=\frac{1}{6}det\left(
\begin{array}{ccc}
x_i+\varepsilon\xi_i&y_i+\varepsilon\eta_i&z_i+\varepsilon\zeta_i\\
x_j+\varepsilon\xi_j&y_j+\varepsilon\eta_j&z_j+\varepsilon\zeta_j\\
x_k+\varepsilon\xi_k&y_k+\varepsilon\eta_k&z_k+\varepsilon\zeta_k
\end{array}
\right)
 $$
 By calculating these determinants for all faces for the total volume we have a presentation:
  $$
 V(P_{\varepsilon})=V_0+\varepsilon V_1+\varepsilon^2 V_2+\varepsilon^3 V_3.
 $$
 By the supposition $V'_{\varepsilon}(0)=0$, so $V_1=0$. If $V_3\ne 0$ then two isometric polyhedra $P_{\varepsilon}$
 and $P_{-\varepsilon}$ have different volumes that is the polynomial (\ref{5}) has two different roots tending under
 $\varepsilon\to 0$  to the same root $V_0$ of the initial volume polynomial. Thus the volume of our $P$ is a multiple root
 of its volume polynomial.

 Now  we consider the case $V_1=V_3=0$ so $V(P_{\varepsilon})$ is
 $$
 V(P_{\varepsilon})=V_0+\varepsilon^2 V_2.
 $$
 We suppose also that $V_0\ne 0$ (because the root $V=0$ is multiple already). Let $Q'_V(0)\ne 0$. Then the equation
 (\ref{5}) determinates $V$ as an analytical implicit function $V=V(l)$ of $|E|$
 variables $\tilde l_{ij}^2=l_{ij}^2+\varepsilon^2 L_{ij}^2$ as independent arguments in some
 \emph{full} neighborhood of values of edge lengths of $P$ which are not related with coordinates of vertices (recall that
 in general only some collections of nonnegative numbers can be presented as squares of lengths of a polyhedron).

  Accordingly \cite{S2} any small diagonal $d$ satisfies an polynomial equation of the form
   \begin{equation}
 \label{6}
 D(l,V,d)= A_0(l,V)d^{2K}+A_1(l,V)d^{2K-2}+...+A_K(l,V)=0,
 \end{equation}
 where coefficients $A_i, 1\le i\le K$ are some polynomials too in squares of edge lengths and square of the polyhedron's
 volume which not all are identically zero. For polyhedra $P_{\varepsilon}$ the coefficients $A_i$ in (ref{6}) are
 presentable as  follows
 $$
 A_i=a_{i0}+a{i1}\varepsilon^2+...+a_{in_i}\varepsilon^{2n_i},
  $$
 and not all coefficients $a_{ij}$ are zero.

\textbf{Lemma 2.}\emph{Any small diagonal of polyhedra $P_{\varepsilon}$  is represented in the form
$d=d_0+ o(\varepsilon), d_0\ne 0$.}

In the proof of lemma one should to distinguish two cases 1) there is at least one coefficient $a_{i0}\ne 0$
and 2) all coefficients $a_{i0}=0$. In the first case we consider for $D$ from (\ref{6}) its derivative
$D'_d(l,V,d)$. If the derivative is not zero then $d$ is expressed from (\ref{6}) as an implicit function
and its Taylor expansion consists only of powers of $\varepsilon^{2m}$. If this derivative is zero then
$d$ satisfies a similar equation of 2 powers less (after the concellation by $d\ne 0$) and one can continue
the same considerations and finally we arrive either to a case with the possibility to present $d$ by a Taylor expansion
with even powers of $\varepsilon$ or to a biquadratic equation. In the case 2) we should reduce all the coefficients
by the maximal common degree $\varepsilon^{2m}$ and we arrive to the case 1).

Now we note that by lemma 1 there exists at least one small diagonal of the form $d=d_0+a\varepsilon, a\ne 0$,
which is in contradiction with lemma 2. So the supposition $Q'_V\ne 0$ is not true.

Let's remark that a seeming theorem should be true for inf. bendable polyhedra in any space $R^n, n>3$ because by \cite{G}
for them there exist volume polynomials too, but for the moment we don't have a needed affirmation about the existence  of
equations for small diagonals similar to (\ref{6}). It would be interesting also to find an algebraic and geometrical
interpretation for the miltiplicity order of the volume as a root of a volume equation.

\end{document}